\documentclass[reqno]{amsart}

\usepackage[all]{xy}

\numberwithin{equation}{section}

\theoremstyle{plain}
\newtheorem{lemma}{Lemma}[section]
\newtheorem{teo}[lemma]{Theorem}

\newtheorem{con}[lemma]{Conjecture}

\theoremstyle{definition}

\theoremstyle{remark}

\newcommand{\pic}{{\rm Pic\thinspace}}
\newcommand{\bs}{{\rm Bs\thinspace}}
\newcommand{\p}{\mathbb{P}}

\newcommand{\st}{\mathbb{S}}

\newcommand{\oc}{{\mathcal O}}
\newcommand{\ls}{{\mathcal L}}
\newcommand{\ms}{{\mathcal M}}

\newcommand{\ci}{{\mathcal I}}
\newcommand{\fb}{\mathfrak{b}}
\newcommand{\cx}{{\mathcal X}}

\newcommand{\rt}{\rightarrow}

\begin{document}

\title{Recent results on linear systems on generic $K3$ surfaces}
\author{Cindy De Volder}
\address{
Department of Pure Mathematics and Computeralgebra,
Galglaan 2, \newline B-9000 Ghent, Belgium}
\email{cdv@cage.ugent.be}
\thanks{The first author is a Postdoctoral Fellow of
the Fund for Scientific Research-Flanders (Belgium) (F.W.O.-Vlaanderen)}

\author{Antonio Laface}
\address{
Dipartimento di Matematica, Universit\`a degli Studi di Milano,
Via Saldini 50, \newline 20133 Milano, Italy
}
\email{antonio.laface@unimi.it}
\thanks{The second author would like to thank the European Research and
Training Network EAGER for the support provided at Ghent University.
He also acknowledges the support of the MIUR of the Italian Government
in the framework of the National Research Project ``Geometry in Algebraic
Varieties'' (Cofin 2002)}
\begin{abstract}
In this note we relate about the problem of evaluate the dimension
of linear systems through fat points defined on generic $K3$
surfaces.
\end{abstract}
\maketitle

\section{Introduction and statement of the problem}

In what follows we assume that the ground field is algebraically
closed of characteristic 0. With $S$ we always denote a smooth
projective {\em generic} $K3$ surface, i.e. $\pic(S)=\langle
H\rangle$ and let $n=H^2$. Consider $r$ points in general position
on $S$, to each one of them associate a natural number $m_i$
called the {\em multiplicity} of the point. We will denote by $\ls
= \ls^n(d,m_1,\ldots, m_n)$ the linear system $|dH|$ through the
$r$ points with the given multiplicities. Define the {\em virtual
dimension} of the system as $v(\ls) = d^2n/2 + 1 -\sum
m_i(m_i+1)/2$ and its {\em expected dimension} by
$e=\max\{v,-1\}$. Observe that $e\leq\dim(\ls)$ and that the
inequality may be strict if the conditions imposed by the points
are dependent. In this case we say that the system is {\em
special}. By $S'$ we will denote the blow-up of $S$ along the $r$
points, given two curves $A,B$ on $S$, the intersection $AB$ will
be defined as the intersection of their strict transforms on $S'$.
The problem of classifying special systems has been largely
studied for linear systems on the plane~\cite{CM2,AG,BS} and more
generally for systems on rational surfaces~\cite{BH,AH}. The main
conjecture on the structure of such systems has been formulated
in~\cite{AH}. In this note we report about some recent results in
the case of generic $K3$ surfaces. In~\cite{CM1} the authors
proved that on the projective plane this conjecture is equivalent
to an older one given by Segre in~\cite{BS}. The advantage of
Segre conjecture is that it can be formulated in the same way on
any surface. Starting from this idea we proved in~\cite{dl1} the
equivalence of Conjecture~\ref{segre} with Conjecture~\ref{our} on
a generic $K3$ surface. An attempt to prove Conjecture~\ref{our}
has been done in~\cite{dl2} by using a degeneration technique
inspired by~\cite{BZ}. The main result, by using this technique,
is Theorem~\ref{deg} which relates the speciality of some linear
systems through points of the same multiplicity with the
speciality of systems through just one point.

\section{The equivalence of the two conjectures}

As stated in the introduction we consider here an extension, to
any surface, of Segre conjecture about special linear systems.

\begin{con}\label{segre}
If $\ls$ is non-empty and reduced linear system on a surface $S$,
then it is non-special.
\end{con}

By Bertini second theorem, this conjecture tell us that if $\ls$
is special, then there exists an irreducible curve $C$ such that
$2C\subseteq\bs(\ls)$. This means that, if Conjecture~\ref{segre}
is true, then in order to give a classification of special systems
on a surface we should be able to classify the type of the curve
$C$. In the case of generic $K3$ surfaces we proved the
equivalence of the preceding conjecture with the following
(see~\cite{dl1}).

\begin{con}\label{our}
Let $\ls$ and $S$ be as above, then
\begin{itemize}
\item[(i)]
    $\ls$ is special if and only if $\ls = \ls^4(d,2d)$ or $\ls = \ls^2(d,d^2)$ with $d \geq 2$;
\item[(ii)]
    if $\ls$ is non-empty then its general divisor has
    exactly the imposed multiplicities in the points $p_i$;
\item[(iii)]
    if $\ls$ is non-special and has a fixed irreducible component $C$
    then
    \begin{itemize}
    \item[a)] $\ls=\ls^2(m+1,m+1,m)=mC+\ls^2(1,1)$ with $C=\ls^2(1,1^2)$ or
    \item[b)] $\ls = 2C$ with $C \in \{ \ls^4(1,1^3),\allowbreak \ls^6(1,2,1),\allowbreak \ls^{10}(1,3)\}$ or
    \item[c)] $\ls = C$.
    \end{itemize}
\item[(iv)]
    if $\ls$ has no fixed components then either its general element is irreducible or
    $\ls=\ls^2(2,2)$.
\end{itemize}
\end{con}

The proof of this result proceeds by analyzing the base locus of
the system $\ls$. Assume that there exist distinct irreducible
curves $C_i$ and $D_j$ such that
\[
\ls = \sum_{i=1}^a\mu_i C_i + \sum_{i=1}^bD_i + \ms,
\]
where $\mu_i\geq 2$ and $\ms$ has no fixed components. By putting
$A,B$ to be two of the irreducible curves into the fixed part of
$\ls$ and assuming conjecture~\ref{segre} to be true, we have that
$v(A)=v(B)=v(A+B)=0$. Since $v(A+B)=v(A)+v(B)+AB-1$, this implies
that $AB=1$. Hence this gives that $C_iC_j=C_iD_j=D_iD_j=1$ and
$C_i^2\leq 1$. Now, it is possible to prove (see~\cite{dl1}) that
given two distinct irreducible curves $A$ and $B$ on $S$ then
either $AB\neq 1$ or $A=\ls^2(1,1^2)$ and $B$ is an irreducible
element of $\ls^2(1,1)$.

\section{A degeneration of K3 surfaces}

In this section we consider an attempt to prove
conjecture~\ref{our} by using a degeneration of $K3$ surfaces to a
union of planes and the blow-up of a $K3$ along points. Let
$\Delta$ be an open disk and let $X$ be the blow-up of
$S\times\Delta$ along $b$ general points of $S\times\{0\}$. The
threefold $X$ is equipped with two projections $p_1,p_2$ on
$\Delta$ and $S$ respectively and the general fiber $X_t$ of $p_1$
is isomorphic to $S$, while $X_0$ is a reducible surface given by
the union of $b$ planes with a surface $\st$. The last surface is
the blow-up of $S$ along the $b$ points. Each one of the $b$
planes $\p_i$ cuts a curve $R_i$ on $\st$ which is a line in
$\p_i$ and a $(-1)$-curve in $\st$. Now given a line bundle $L$ on
$S$ it is possible to construct infinitely many line bundles
(depending on the integer $k$) $ \oc_{X}(L,k) := p_2^* (L) \otimes
\oc_X(k\st)$ on $X$ such that each one restricted to $X_t$ gives
$L$. Defining $\cx(L,k)$ as the restriction to $X_0$ we have that
\[
\begin{array}{lcl}
\cx(L,k)_{|\p^i} & = & \oc_{\p^2}(k) \\
\cx(L,k)_{|\st} & = & \fb^*(L) \otimes \oc_{\st}(- \sum_{i=1}^{b}
k E_i),
\end{array}
\]
where $\fb:\st\rt S$ is the blow-up map. This construction allows
us to degenerate a system on $S$ to a union of systems on the
$\p_i$'s and $S$ in the following way. Let $Z := m_1q_1 + \cdots +
m_r q_r$ be a subscheme of $S$ with points in general position.
Chosen $a_1, \ldots, a_b$ positive integers such that $a_1 +
\cdots + a_b \leq r$, let $Z'_i$ be the specialization of $a_i$
points of $Z$ to points of $\p_i$ (with the same multiplicities).
Let $Z'_{\st}$ be the residual subscheme, made of  $r - \sum a_i$
general points of $\st$. Given $Z' := Z'_1+\ldots+Z'_b+Z'_{\st}$,
one has that $\cx(\ls,k)\otimes \ci_{Z'}$, is a degeneration of
$\ls\otimes\ci_Z$. In this way, the starting system $\ls$ through
$r$ degenerate to the system $\ls_0$ on $X_0$ made by the $\ls^i$
on the $\p_i$ and by the $\ls_{\st}$ on $\st$. Observe that the
last system corresponds to a system on $S$ through less than $r$
points. In this way, by using the fact that the homogeneous planar
systems $\ls_2(d,m^4),\ \ls_2(d,m^9)$ are never special, it is
possible to use the preceding degeneration in an inductive way.
So, for example consider the system $\ls^n(d,m^{4^h})$, take
$b=4^{h-1}$ and put four general points on each of the $\p_i$. In
this way the speciality of the starting system is related to that
of $\ls^n(d,m^{4^{h-1}})$ and so on. More generally we have the
following (see~\cite{dl2}).

\begin{teo}\label{deg}
If $\ls^{n}(d,m)$ is non-special for all non-negative integers
$(d,m)$ then $\ls^{n}(d',m'^{4^h 9^k})$ is non-special for all
non-negative integers $(d',m',h,k)$.
\end{teo}

Unfortunately it is an open problem to evaluate if a system
through just one point is special or not. The only known example
is $\ls^4(d,2d)$ as stated in Conjecture~\ref{our}.

\newpage

\bibliographystyle{alpha}


\end{document}